\documentclass[transmag]{IEEEtran}

\usepackage{cite}
\usepackage{mathtools,amssymb,amsfonts}		

\usepackage{graphicx}
\graphicspath{{./figs/}}
\DeclareGraphicsExtensions{.pdf,.png,.jpeg,.jpg}

\def\BibTeX{{\rm B\kern-.05em{\sc i\kern-.025em b}\kern-.08em
    T\kern-.1667em\lower.7ex\hbox{E}\kern-.125emX}}

\usepackage[hyphens]{url}
\usepackage[T1]{fontenc}
\usepackage[utf8]{inputenc}
\usepackage{siunitx}
\DeclareSIUnit{\sample}{S}
\DeclareSIUnit{\belmilliwatt}{Bm}
\DeclareSIUnit{\dBm}{\deci\belmilliwatt}

\usepackage{xcolor}
\usepackage{acronym}

\usepackage{algorithm, algpseudocode}
\usepackage{multirow}

\usepackage{subcaption}

\usepackage{bm}
\usepackage{upgreek}

\DeclareRobustCommand{\mathbfup}[1]{\begingroup\changegreekbf\mathbf{#1}\endgroup}

\makeatletter
\def\changegreek{\@for\next:={%
		alpha,beta,gamma,delta,epsilon,zeta,eta,theta,kappa,lambda,mu,nu,xi,pi,rho,sigma,%
		tau,upsilon,phi,chi,psi,omega,varepsilon,vartheta,varpi,varrho,varsigma,varphi}%
	\do{\expandafter\let\csname\next\expandafter\endcsname\csname up\next\endcsname}}
\def\changegreekbf{\@for\next:={%
		alpha,beta,gamma,delta,epsilon,zeta,eta,theta,kappa,lambda,mu,nu,xi,pi,rho,sigma,%
		tau,upsilon,phi,chi,psi,omega,varepsilon,vartheta,varpi,varrho,varsigma,varphi}%
	\do{\expandafter\def\csname\next\expandafter\endcsname\expandafter{%
			\expandafter\bm\expandafter{\csname up\next\endcsname}}}}
\makeatother

\newcommand*{\figref}[1]{Fig.~\ref{#1}}
\newcommand*{\tabref}[1]{Table~\ref{#1}}
\newcommand*{\secref}[1]{Section~\ref{#1}}

\newcommand*{\abs}[1]{\left|#1\right|}
\newcommand*{\norm}[1]{\left\lVert#1\right\rVert}

\newcommand*{\vect}[1]{\ensuremath{\mathbfup{#1}}}
\newcommand*{\matx}[1]{\ensuremath{\mathbfup{#1}}}

\newcommand{\revision}[1]{#1}

\acrodef{kaf}[KAF]{kernel adaptive filtering}
\acrodef{klms}[KLMS]{kernel least mean squares}
\acrodef{krls}[KRLS]{kernel recurusive least squares}
\acrodef{kapa}[KAPA]{kernel affine projection algorithm}
\acrodef{lms}[LMS]{least mean squares}
\acrodef{rls}[RLS]{recursive least squares}
\acrodef{ald}[ALD]{approximate linear dependency}
\acrodef{ls}[LS]{least squares}
\acrodef{im2lms}[IM2LMS]{second-order intermodulation least mean squares}

\acrodef{sgd}[SGD]{stochastic gradient decent}
\acrodef{rx}[Rx]{receive}
\acrodef{tx}[Tx]{transmit}
\acrodef{rf}[RF]{radio frequency}
\acrodef{dac}[DAC]{digital-to-analog converter}
\acrodef{adc}[ADC]{analog-to-digital converter}
\acrodef{pa}[PA]{power amplifier}
\acrodef{snr}[SNR]{signal-to-noise ratio}
\acrodef{lna}[LNA]{low noise amplifier}
\acrodef{lo}[LO]{local oscillator}
\acrodef{fdd}[FDD]{frequency division duplex}
\acrodef{csf}[CSF]{channel select filter}
\acrodef{ca}[CA]{carrier aggregation}
\acrodef{osf}[OSF]{oversampling factor}
\acrodef{fir}[FIR]{finite impulse response}
\acrodef{dc}[DC]{direct-current}
\acrodef{imd2}[IMD2]{second-order intermodulation distortion}
\acrodef{rb}[RB]{resource block}
\acrodef{i}[I]{in-phase}
\acrodef{q}[Q]{quadrature}
\acrodef{lte}[LTE]{long term evolution}

\acrodef{tdl}[TDL]{tapped delay line}
\acrodef{lut}[LUT]{look-up table}

\acrodef{nmse}[NMSE]{normalized mean squared error}
\acrodef{snr}[SNR]{signal-to-noise ratio}
\acrodef{sinr}[SINR]{signal-to-interference-and-noise ratio}
\acrodef{icn}[ICN]{interference-to-noise ratio}

\acrodef{svm}[SVM]{support vector machine}

\acrodef{cpd}[CPD]{canonical polyadic decomposition}

\begin{document}
	
	\title{Enhanced Nonlinear System Identification by Interpolating Low-Rank Tensors}
	
	\author{Christina Auer, \IEEEmembership{Graduate Student Member, IEEE}, 
		Thomas Paireder, \IEEEmembership{Graduate Student Member, IEEE},\\ 
		Oliver Ploder, \IEEEmembership{Graduate Student Member, IEEE}, Oliver Lang,
		\IEEEmembership{Member, IEEE}, \\
		and Mario Huemer,
		\IEEEmembership{Senior Member, IEEE}
		\thanks{Acknowledgments: The financial support by the Austrian Federal Ministry for Digital and Economic
			Affairs, the National Foundation for Research, Technology and Development and the Christian Doppler
			Research Association is gratefully acknowledged.}
		\thanks{C. Auer, T. Paireder and O. Ploder are with the Christian Doppler Laboratory for Digitally Assisted RF Transceivers for Future Mobile Communications, Institute of Signal Processing, Johannes Kepler University, Linz, Austria (e-mail: christina.auer@jku.at).}
		\thanks{O. Lang and M. Huemer are with the Institute of Signal Processing, Johannes Kepler University, Linz, Austria}}

	\IEEEtitleabstractindextext{\begin{abstract} 
			Function approximation from input and output data is one of the most investigated problems in signal processing. This problem has been tackled with various signal processing and machine learning methods. Although tensors have a rich history upon numerous disciplines, tensor-based estimation has recently become of particular interest in system identification. In this paper we focus on the problem of adaptive nonlinear system identification solved with interpolated tensor methods. We introduce three novel approaches where we combine the existing tensor-based estimation techniques with multidimensional linear interpolation. To keep the reduced complexity, we stick to the concept where the algorithms employ a Wiener or Hammerstein structure and the tensors are combined with the well-known \ac{lms} algorithm. The update of the tensor is based on a stochastic gradient decent concept. Moreover, an appropriate step size normalization for the update of the tensors and the \ac{lms} supports the convergence. Finally, in several experiments we show that the proposed algorithms almost always clearly outperform the state-of-the-art methods with lower or comparable complexity.  
		\end{abstract}
		
	}

\maketitle


\section{INTRODUCTION}

\IEEEPARstart{I}{n} recent decades, great attention in the scientific community has been directed to the problem of system identification. For this problem machine learning-based methods gained tremendous popularity. Machine learning is mainly associated with deep learning and neural networks ~\cite{liu2017survey,5406124,ploderVTC2019,BalatsoukasStimming2018NonLinearDS,ploderAsilomar2019}, but this field covers many other techniques such as random forests~\cite{hastie01statisticallearning,8827054}, support vector machines~\cite{auerVTC2020,Smola2004,cortes_support-vector_1995}, kernel adaptive filters~\cite{txhamchristina,liu_principe_haykin_kernel_adaptive_filtering, Kennedy-book:2013} and tensor-based learning (cf. e.g. \cite{cichocki_tensor_2015}).

Although tensor-based methods can perform similarly to other mentioned techniques~\cite{sidiropoulos_tensor_2017,kargas_nonlinear_2019} they have been ignored for a long time, possible due to their high computational complexity. However, tensor-based methods have recently gained more and more attention. Tensor-based algorithms have been used for many tasks such as face recognition, mining musical scores and detecting cliques in social networks~\cite{10.1007/3-540-47969-4_30,cichocki_tensor_2015,10.1007/978-3-642-33460-3_39}.  Moreover, they have been already deployed in system identification~\cite{bousse2017tensor,favier2008tensor,6190040}, channel identification~\cite{fernandes2009blind,kibangou2007blind}, separation of speech signals~\cite{5233821}, emitter localization for radar, communication applications and passive sensing ~\cite{5510180,852018}. 

The issue of complexity has been addressed in various forms and papers. Tensor trains~\cite{ko_fast_2015} were developed as an attempt to make the overall system less complex and less memory intensive. In order to reduce the necessary dimension of the tensor, it is split into multiple smaller tensors, cascading one another. In \cite{batselier_extended_2019, batselier_matrix_2017}, an interesting combination of a tensor and an extended Kalman filter was investigated with the goal of enhancing the performance and keeping the rank of the tensor low. Further, a combination of tensor trains and B-splines was used in~\cite{karagoz_nonlinear_2020}   to mitigate the problem of complexity, while also leveraging the advantages of B-splines for their system identification task. 

In~\cite{tensor}, which is also the basis of this paper, the approach of Scarpiniti et.al.~\cite{splinepaperelsevier} was used for system identification by splitting  up the systems following a Wiener or Hammerstein model~\cite{wiener_1958} into their linear and nonlinear parts. However, a big downside is, that these problems are continuous. The tensor is basically a huge \ac{lut} and only a few discretization points, which define the size of the tensor, might not be enough. Hence, the output might not always be sufficiently accurate and thus, the performance can suffer a lot. 

In this work we are also investigating system identification and we assume that the unknown system changes over time. We make use of the adaptive tensor combination approaches in \cite{tensor}. The way we mitigate the problem of inaccuracy is to combine these classical tensor methods with interpolation. 
In fact, we introduce three novel interpolated tensors methods, the interpolated tensor-only, a Wiener and a Hammerstein system. The latter two combine the tensor variant with \ac{fir} filters, which are trained with the well-known \ac{lms} algorithm and which model the linear part of the system. Interpolation is achieved the following way. First, the elements of the tensor, whose corresponding input samples are closest to the input sample for which the interpolation should be \revision{carried} out, are stored in a separate subtensor of dimension $2 \times 2 \times \hdots \times 2$. Second, the elements of this subtensor are interpreted as control points, and they are linearly interpolated to get a more accurate output. The \ac{sgd} method is used for the tensor update. Moreover, convergence is ensured by an appropriate novel step size normalization. To conclude the investigations numerical experiments have been performed which verified the huge increase in accuracy provided by the interpolation.
\revision{
In summary the novel contributions of this paper are:
\begin{itemize}
	\item the mathematical framework for interpolated tensors
	\item three models using the interpolated tensors methods:
	\begin{itemize}
		\item the interpolated decomposed tensor model
		 \item the interpolated tensor-LMS model
		 \item the interpolated LMS-tensor model and
		\end{itemize}
	\item the step size normalization.
\end{itemize} }

This paper is organized as follows. For completeness, in \secref{sec:basics}, we present some essential notation and basics needed in the remainder of this paper. In \secref{sec:sys_ident}, we derive the novel interpolated tensor models, and in \secref{sec:step_size} we introduce a constraint to ensure convergence. Then, \secref{complexity} describes the computational complexity while \secref{sec:results} shows some experimental results. Finally, \secref{conclusion} concludes this work. 

\section{Basics}

\label{sec:basics}

In physics and mathematics the word tensors can refer to different meaning. In this paper, we mean by a tensor a multilinear mapping between two sets as widely adopted and also used in \cite{sidiropoulos_tensor_2017,kargas_nonlinear_2019}. As a basis, we introduce the notation adopted from \cite{tensor} and recapitulate the most important properties and definitions of the well-known tensor methods from \cite{tensor}. Moreover, we give a very brief overview of multidimensional linear interpolation. These basics are then utilized in the derivations in \secref{sec:sys_ident}.

\subsection{Preliminaries and Notation}
\revision{
A tensor $\mathcal X \in \mathbb R^{I_1 \times \dots\times I_M} $ may be represented as an $M$-dimensional array, indexed by $\mathcal X( i_1,i_2,i_3, \dots, i_M)$. 
Multiplying a tensor  $\mathcal X \in \mathbb R^{I_1 \times \dots\times I_M}$  by a vector $\mathbf v  \in \mathbb R^{ I_m \times 1}$ can be done via the $m$-mode product for vectors
\begin{align}
&\left (\mathcal X \times_m \mathbf v^{\text{T}} \right)(i_1,\dots,i_{m-1},1,i_{m+1}, \dots, i_M) \nonumber \\
& =\sum_{i_m = 1}^{I_m}  \mathcal X (i_1,\dots,i_M) \mathbf v(i_m).
\end{align}
Note that the $m$-th dimension of the resulting tensor is 1, which means that the order of the tensor has reduced to $M-1$.
}
\revision{
The canonical polyadic decomposition (CPD) 
\begin{align}
&\mathcal X (i_1,i_2,\dots,i_M) =  \sum_{r = 1}^R \prod_{m = 1}^M \mathbf A_m(i_m,r),
\label{eq:forwardtensor}
\end{align}
where $R$ denotes the rank,
can be seen as a generalization of the matrix singular value decomposition. 
}

\revision{
Reshaping a tensor into a matrix is called matricization, unfolding or flattening. This can be done along all coordinates.  Let $m'$ be an arbitrary but fixed number $m' \in \{1, \dots, M\}$. Then the resulting matrix is denoted by $\mathcal X_{(m')}$, where the index $(m')$ indicates how the tensor is flattened. A cube for example in the three-dimensional case can be cut horizontally, in front slabs, or sideways. Regardless of the dimension the resulting matrices lead to reshaping the tensor into one matrix $\mathcal X_{(m')} \in \mathbb R^{I_1 \cdots I_{m'-1}I_{m'+1} \cdots I_M \times I_{m'}}$. Introducing some common notation as in \cite{sidiropoulos_tensor_2017,kargas_nonlinear_2019}, 
we denote the Hadamard product by $\circledast $ and the product 
over all matrices $\mathbf A_{m}$ with $m \neq m'$ as
\begin{align}
\label{eq:hadamard}
\circledast_{m\neq m'} \mathbf A_{m}: = \mathbf A_M \circledast \cdots \circledast \mathbf A_{m'+1} \circledast \mathbf A_{m'-1} \circledast \cdots \circledast \mathbf A_1 .
\end{align}
}

\subsection{Tensor Decomposition}

With the decomposed tensor, as described above, we can investigate the \ac{sgd} algorithm \cite{kargas_nonlinear_2019} by minimizing the following cost function

\begin{align}
J_{\text{SGD}} = \left(\mathcal Y(i_1,i_2,\dots,i_M) - \sum_{r = 1}^R \prod_{m = 1}^M \mathbf A_m(i_m,r)\right)^2 =: e^2
\end{align}
where $\mathcal Y$ is a tensor containing the corresponding data points. For solving this minimization problem, we need to calculate the gradient

\begin{align}
\frac{\partial J_{\text{SGD}}}{\partial \mathbf A_{m'}(i_{m'},:)  } 
= - 2e \circledast_{m\neq m'} \mathbf A_m(i_m,:) .
\end{align}
and can therefore, formulate the following update equation 
\begin{multline}
\label{eq:SGD_update}
\mathbf A_{m',n+1}(i_{m',n+1},:) = \\
\hspace{3em}\mathbf A_{m',n}(i_{m',n},:) -\mu \left(\frac{\partial J_{\text{SGD}}}{\partial \mathbf A_{m'}(i_{m'},:)}\right) ^{\text{T}}  
\end{multline}
where $\mu$ is the step size. From the available data points, we randomly take one $\mathcal Y(i_1,i_2,\dots,i_M)$ and only take the gradient step for those model parameters that have an effect. 
\revision{
Since we are talking about an update step an additional iteration parameter $n$ is necessary. For clarity it was skipped above \eqref{eq:SGD_update} and will be skipped whenever it is not necessary. Later in the problem description $n$ will also be the time index for the appearing time signals. 
}

\subsection{Multidimensional Linear Interpolation}

\label{sec:multdim_lin_intpol}
Linear interpolation is a method to construct new data points within the range of a discrete set of known data points using linear polynomials. This can be seen as finding a continuous function with this discrete data set. We call multidimensional linear interpolation an extension for using repeated linear interpolation of functions with various variables on a regular grid. It approximates the value of a function $f$ at an intermediate point $(\check v_1, \check v_2,\dots \check v_M)$  within that grid. 
The function $f$ to be interpolated is known at given points $(v_{1,n},v_{2,n},\dots v_{M,n})$  and the goal is to find the values of the function at arbitrary points $(\check v_1, \check v_2,\dots \check v_M)$. 
Considering the bilinear interpolation (= interpolation of two variables using repeated linear interpolation)  on a unit square this would be 
\begin{align}
f(\check v_1, \check v_2) \approx
\begin{bmatrix}
1- \check v_1 & \check v_1 
\end{bmatrix}
\begin{bmatrix}
f(0,0) & f(0,1) \\
f(1,0) &f(1,1)
\end{bmatrix}
\begin{bmatrix}
1- \check v_2\\ \check v_2
\end{bmatrix}.
\end{align}
For a multidimensional function with $M$ dimensions this would expand to 
\begin{align}
f(\check v_1, \ldots, \check v_M) \approx \sum_{l_1, \ldots, l_M = 0}^{1} f(l_1\dots l_M)
\prod_{m = 1}^{M}
\check {\mathbf v}_{m} (l_m),
\end{align}
where 
$\check{\mathbf v}_m = \begin{bmatrix}1-\check v_m\\\check v_m
\end{bmatrix}$.
Note that the sum indicates a multiple sum, therefore $M$ ones over each $l_m \in \{0,1\}$.

\subsection{Problem Statement}
\begin{figure}[t!]
	\centering
	\includegraphics[width=0.75\columnwidth]{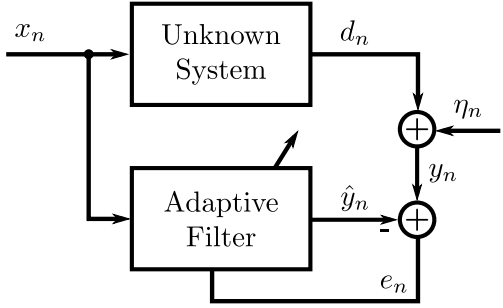}
	\caption{A general setup of adaptive system identification. \revision{The unknown system is estimated by the adaptive filter. In this paper it follows a Wiener or a Hammerstein structure.} Figure taken from \cite{tensor}.}
	\label{fig:overviewarch}
\end{figure}

The aim of this paper is to approximate an unknown system, 
which can be seen in the general setup in Figure~\ref{fig:overviewarch}. 
Therefore, we use the same input signal $x_n$, \revision{ where $n$ indicates the time index, } as the unknown system and observe the noisy output $y_n = d_n+\eta_n$, where $d_n$ is the desired signal, $\eta_n$ is the noise. The goal is to minimize the error $e_n = y_n -\hat{y}_n$, where $\hat{y}$ is the output of the adaptive filter and $\hat{(\cdot)}$ always defines the approximation. As a consequence of the additive noise, this error can never truly equal zero.  
Since we assume that the unknown system may not remain static over time, the adaptation will not be turned off and the adaptive filter is updated with each observed sample (i.e. one adaption step per time-step).  
Note, that the unknown system we want to approximate can also have a time memory and therefore, not only depends on the actual input signal $x_n$ but also on its past samples. 
Hence, the $P$ most recent samples of $x_n$ are stored in a \ac{tdl} according to
\begin{align}
\mathbf x_n: = \text{tdl}(x_n,\mathbf x_{n-1} ) &:= [x_n,\mathbf x_{n-1}^{\text T}(1:P{-}1)]^{\text T} 
\end{align}
with $\mathbf x_0 = \mathbf 0$ and $\mathbf x_n\in \mathbb R^{P}$.

\begin{figure}[t!]
	\centering
	\begin{subfigure}[t]{\textwidth}
		\includegraphics[width=0.35\textwidth]{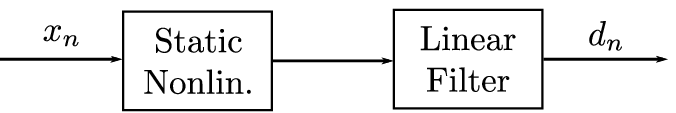}
		\caption{}
		\label{fig:ham}
	\end{subfigure}
	~ 
	\begin{subfigure}[t]{\textwidth}
		\includegraphics[width=0.35\textwidth]{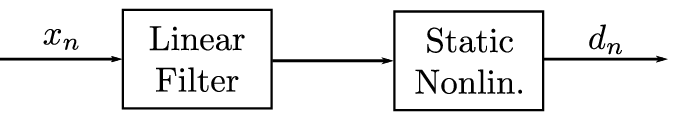}
		\caption{}
		\label{fig:wiener}
	\end{subfigure}
	\caption{\revision{Considered unknown system following either (a) a Hammerstein, or (b) a Wiener structure. Figure taken from \cite{tensor}.} }
	\label{fig:lin_nonlin}
\end{figure}

\revision{Later in this paper, we will investigate six different systems which we will approximate with the different methods. Three of these unknown systems, cf. \figref{fig:overviewarch}, are assumed to following a Hammerstein structure, cf. \figref{fig:ham}, while the other three systems area assumed to follow a Wiener structure, cf. \figref{fig:wiener}}.

\subsection{Discretization}

The input of the decomposed tensors are the indices of the matrices $\mathbf A_m$.  
Hence, naturally there can only be $I_m$ (the size of  first dimension of $\mathbf A_m$) different input values and not arbitrary many. 
Since our input signal $x_{n}$ is usually not limited to $I_m$ different values, it has to pass a separate discretization step. 
This is carried out by dividing the range between the maximum and minimum of the input values in a number of intervals, which then also defines the first dimension of the matrix $\mathbf A_{m}$. 
A certain cut-off point can be defined in case the maximum and/or minimum values are unknown.
In case $I_m$ is even, this discretization step is mathematically written as
\begin{align}
\label{eq:disc}
\text{disc}(x_n) &= \left\lfloor\frac{x_n}{\Delta x}\right\rfloor + \frac{I_{m} }{2}
\end{align}
where $\Delta x$ is the discretization interval. $\Delta x$ is kept the same for all $M$ dimensions of the tensor for the sake of simplicity. An extension to different discretization intervals for the different dimensions is straight forward. The floor operator $\lfloor \cdot \rfloor$ is defined the following way
\begin{align} 
\lfloor \cdot \rfloor: = \max\{k\in \mathbb Z \mid k\leq x\}.
\end{align}

\section{Adaptive System Identification with Interpolated Tensors}
\label{sec:sys_ident}

\revision{
We first describe the basic tensor model. Based on that, we introduce three algorithms: 
\begin{itemize}
	\item the interpolated decomposed tensor model 
	\item the tensor-LMS model, and 
\item the LMS-tensor model.
\end{itemize}
}
\begin{figure*}[t!]
	\centering
	\begin{subfigure}[t]{\textwidth}
		\centering
		\includegraphics[width=0.8\textwidth]{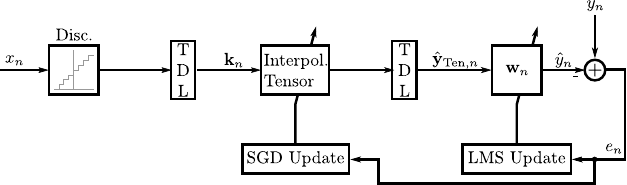}
		\caption{\revision{
				First the input samples pass through the discretization step and the tapped delay line. Then the interpolated tensor steps follows. Since the input of the LMS needs to be a vector again a TDL  precedes this step. 
			}}
		\label{fig:mod1}
	\end{subfigure}
	~ 
	\begin{subfigure}[t]{\textwidth}
		\centering
		\includegraphics[width=0.8\textwidth]{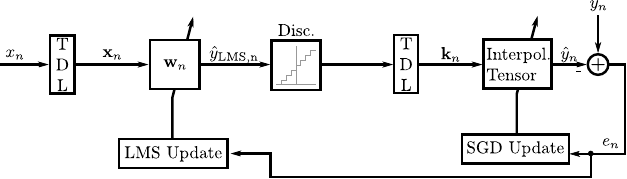}
		\caption{\revision{Here the input samples take similar steps with the LMS and tensor swapped. The inputs samples are first gathered by the TDL followed by the LMS. The output is then discretized, and before passing the interpolated tensor the samples have to go through the TDL.
		}}
		\label{fig:mod2}
	\end{subfigure}
	\caption{The (a) tensor-LMS, and (b) LMS-tensor architecture taken from\cite{tensor} for the Hammerstein and Wiener model, where TDL denotes the tapped delay line.  
	 }
	\label{fig:models}
\end{figure*}

\subsection{Basic Tensor Model}
The goal is to estimate the value of the function $f(x_{1}, \ldots, x_{M})$ by linearly interpolating control points, which are elements of a tensor $\bm{\mathcal{Q}}$.
Each entry of that tensor corresponds to a function value, where the first dimension of the function is $x_{1}$, second is $x_{2}$, etc.
\begin{equation}
\bm{\mathcal{Q}}(k_1, \ldots, k_M) = f\left(\left(k_1 - \frac{I_1}{2}\right) \Delta x, \ldots\right)
\end{equation}
with $k_m  = \text{disc}(x_{m})$
for all $m=1,\dots,M$. Out of this total tensor we consider a subtensor 
\begin{align}
\label{eq:tensor_points}
\bm{\mathcal{Q}}_{k_1, \ldots, k_M} &\in \mathbb{R}^{2 \times 2 \times \ldots \times 2}
\end{align}
which consists of the point described in \eqref{eq:tensor_points} and its neighboring ones, for each function value to be estimated. 

Using this knowledge and the mathematical basics about multidimensional linear interpolation described in \secref{sec:multdim_lin_intpol} we can now estimate any value 
of the function $f$ by generally interpolating the values from the subtensor
\begin{align}
&f(x_1, \ldots, x_M) \approx \bm{\mathcal{Q}}_{k_1, \ldots, k_M} \times_M \check{\vect{u}}_M \ldots \times_1 \check{\vect{u}}_1 = \nonumber \\
&  \sum_{l_1, \ldots, l_M = 0}^{1} \bm{\mathcal{Q}}_{k_1, \ldots, k_M}(l_1, \ldots, l_M) \prod_{m = 1}^{M} \check{\vect{u}}_m(l_m) =: \hat y_{\text{SGD}}
\label{eq:intpol}
\end{align}
with 
\begin{align}
\check{\vect{u}}_m &=
\begin{bmatrix}
1 - u_m\\
u_m
\end{bmatrix}\\
u_m &= \frac{x_m}{\Delta x} - \left\lfloor\frac{x_m}{\Delta x}\right\rfloor.
\end{align}


One approach for updating the subtensor $\bm{\mathcal{Q}}_{k_1, \ldots, k_M}$, and with that also the full tensor $\bm{\mathcal{Q}}$, is the \ac{sgd} approach. Here, we utilize the squared error cost function
\begin{align}
J & = \left(y_n - \hat y_{\text{SGD}}\right)^2 := e_n^2.
\end{align}
In order to solve the minimization problem, we calculate the derivative
\begin{align}
&\frac{\partial J}{\bm{\mathcal{Q}}_{k_{1,n}, \ldots, k_{M,n},n}(l'_1, \ldots, l'_M)} = 
- 2 \, e_n \, \prod_{m = 1}^{M} \check{\vect{u}}_{m,n}(l'_m) \\
& := - 2 \, e_n \, \bm{\mathcal{B}}_{k_{1,n}, \ldots, k_{M,n},n}(l'_1, \ldots, l'_M), 
\nonumber
\end{align}
for each combination of $ l'_m \in \{0,1\}$.
With this derivative, we get the following update formulation of the whole subtensor 
\begin{align}
\bm{\mathcal{Q}}_{k_{1,n}, \ldots, k_{M,n},n+1} = 
\bm{\mathcal{Q}}_{k_{1,n}, \ldots, k_{M,n},n} + \nonumber\\
2\mu  \, e_n \,  \bm{\mathcal{B}}_{k_{1,n}, \ldots, k_{M,n},n}.
\end{align}

\subsection{Interpolated Decomposed Tensor Model}
\label{sec:intpol_decomp}

We recall the regular tensor decomposition \eqref{eq:forwardtensor} which has been studied in several publication \cite{sidiropoulos_tensor_2017,kargas_nonlinear_2019,tensor} 
\begin{equation}
\bm{\mathcal{Q}} (i_1, \ldots, i_M) = \sum_{r = 1}^{R} \prod_{m = 1}^{M} \matx{A}_m (i_m,r).
\label{eq:tensor_decomp}
\end{equation}
Plugging \eqref{eq:tensor_decomp} into \eqref{eq:intpol}, we can approximate function values with multidimensional linear interpolation discussed in \secref{sec:multdim_lin_intpol}
\begin{align}
f(x_1, \ldots, x_M) \approx \sum_{l_1, \ldots, l_M = 0}^{1} \sum_{r = 1}^{R} \prod_{m = 1}^{M} \matx{A}_m(k_m + l_m,r) \check{\vect{u}}_m(l_m).
\label{eq:estimate_intpol}
\end{align}
Furthermore we can use the \ac{sgd} learning to update the matrices $\matx{A}_{m}$. Therefore, we need the gradient of the cost function where the estimated function value is calculated as in \eqref{eq:tensor_decomp}. The gradient can be observed by taking the following partial gradient for all $m'\in [1,\dots,M]$ and $l'_{m'} \in \{0,1\}$
\begin{align}
\label{eq:update_intpol}
&\frac{\partial J}{\partial \matx{A}_{m'}(k_{m'} + l'_{m'}, :)} =  -2 \, e \, \check{\vect{u}}_{m'}(l'_{m'}) \\ 
&\sum_{l_1, \ldots, l_{m'-1}, l_{m'+1}, \ldots, l_{M} = 0}^{1} \circledast_{m \neq m'} \matx{A}_m(k_{m} + l_m, :) \check{\vect{u}}_m(l_m), \nonumber 
\end{align}
\revision{
where $\circledast_{m \neq m'} $ denote the Hadamard product 
over all matrices $\mathbf A_{m}$ with $m \neq m'$ according to
\eqref{eq:hadamard}.
}

\subsection{Tensor-LMS (TLMS)}

In this subsection we assume that the unknown system in \figref{fig:overviewarch} consists of a nonlinearity and a linear filter, following a Hammerstein model.
For estimating this nonlinear system we use the tensor-LMS (TLMS) architecture shown in \revision{\figref{fig:mod1}} to identify the system. Here, the nonlinearity is modeled by a tensor, where the tensor output is interpolated as in \secref{sec:intpol_decomp}. \revision{ The output samples of the tensor are gathered by the TDL before a linear \ac{fir} filter follows,} which models the linear part, and whose filter coefficients are updated by the \ac{lms} algorithm.


The goal is to estimate the value at the output of the linear filter but the input first goes through the interpolated tensor decomposition, cf. \revision{\figref{fig:mod1}}, resulting in two forward path blocks. 
The cost function is the instantaneous squared error (as usual for the \ac{lms} filter) 
\begin{equation}
J_\text{TLMS} = \left(y_n - \hat y_{\text{LMS},n}  \right)^2 =: e_{n}^2
\end{equation}
with
\begin{align}
\hat y_{\text{LMS},n} = \mathbf w_n^{\text T} \hat{\mathbf y}_{\text{Ten},n}  = \sum_{p = 1}^{P} w_{n,p} \hat y_{\text{Ten},n-p+1},
\end{align}
where $\hat{\mathbf y}_{\text{Ten},n} : = \text{tdl}(\hat y_{\text{Ten},n})$ represents the interpolated tensor output.
The interpolated tensor output is given by 
\begin{equation}
\hat y_{\text{Ten},n} = \sum_{l_1, \ldots, l_M = 0}^{1} \sum_{r= 1}^{R} \prod_{m = 1}^{M} \matx{A}_{m,n}(k_{m,n} + l_m,r) \check{\vect{u}}_{m,n}(l_m).
\end{equation}
Therefore, we can reformulate the cost function
\begin{align}
\begin{split}
&J_\text{TLMS} = \left(y_n -  \sum_{p = 1}^{P} w_{n,p} \sum_{l_1, \ldots, l_M = 0}^{1} \sum_{r= 1}^{R} \right. \\
& \left. \prod_{m = 1}^{M} \matx{A}_{m,n-p+1}(k_{m,n-p+1} + l_m,r) \check{\vect{u}}_{m,n-p+1}(l_m)   \right)^2 \\
&=: \left(y_n -s_{n-p+1} \right) = e_{n}^2.
\end{split} 
\end{align}
In order to be able to calculate the update for both the tensor and the weights of the \ac{lms} we need to calculate the derivatives in the backward path. 
We start over with the tensor derivation
\begin{align}
& \frac{ \partial J_\text{TLMS}}{\partial \mathbf A_{m'}} =  -2 \, e_{n} \,    \sum_{p = 0}^P w_p \nonumber \\
&\left (\begin{pmatrix}
\mathbf 0_{k_{m',n-p+1} - 1\times R}\\
\frac{\partial s_{m',n-p+1}}{\partial \mathbf A_{m'}(k_{m',n-p+1} ,:)}\\
\mathbf 0_{I_{m'} - k_{m',n-p+1}  \times R}\\
\end{pmatrix} + 
\begin{pmatrix}
\mathbf 0_{k_{m',n-p+1} \times R}\\
\frac{\partial s_{m',n-p+1}}{\partial \mathbf A_{m'}(k_{m',n-p+1} + 1,:)}\\
\mathbf 0_{I_{m'}-k_{m',n-p+1} -1\times R}\\
\end{pmatrix} 
\right) \nonumber \\
&=: -2 \, e_{n} \, \mathbf S_{m',n}
\label{eq:derivTLMS_intpolA}
\end{align}
where $\mathbf 0_{k_{m',n-p+1}-1\times R} \in \mathbb R^{k_{m',n-p+1}-1 \times R}$ denotes a matrix with all elements being zero and 
with 
\begin{align}
&\frac{\partial s_{m',n-p+1}}{\partial \mathbf A_{m'}(k_{m',n-p+1} + l'_{m'},:)} = \\
&\check{\vect{u}}_{m',n-p+1}(l'_{m'})\sum_{l_1, \ldots, l_{m'-1}, l_{m'+1}, \ldots, l_{M} = 0}^{1} \nonumber   \\
&\circledast_{m\neq m'} \mathbf A_{m,n-p+1}(k_{m,n-p+1} + l_m,:) \check{\vect{u}}_{m,n-p+1}(l_m). \nonumber
\end{align}
The backward path of the \ac{lms} it quite straight forward
\begin{equation}
\frac{\partial J_\text{TLMS}}{\partial \mathbf{w}_n} = -2 \,  e_{n} \, \hat{\mathbf y}_{\text{Ten},n}.
\end{equation}
Note that since  $ \mathbf{w}_n$ is a column vector and we are using the numerator-layout such that  $\frac{\partial J_\text{TLMS}}{\partial \mathbf{w}_n} $ is a row vector.
With these equations we can calculate the update for the decomposed tensor model as 
\begin{align}  
\label{eq:updateTLMS_intpolA}
&\mathbf A_{m',n+1}
=
\mathbf A_{m',n}
+ 2 \mu_2 e_{n} \,  \mathbf S_{m',n}.
\end{align}
The update for the \ac{lms} is commonly known \cite{kay} and given by 
\begin{align}
\mathbf{w}_{n+1} = \mathbf{w}_{n} + 2 \mu_1  \, e_{n} \, \hat{\mathbf y}_{\text{Ten},n}.
\end{align}

\subsection{LMS-Tensor (LMST)}
For dealing with a problem, where the unknown system in \figref{fig:overviewarch} can be represented by a linear term followed by a nonlinearity we propose an interpolated \ac{lms}-tensor (LMST) architecture, cf. \revision{\figref{fig:mod2}}. 
The linear filter is again modeled by an \ac{fir} filter, whose filter coefficients are updated by an \ac{lms} algorithm, while the nonlinearity is modeled by a tensor whose output is interpolated as described in \secref{sec:intpol_decomp}.
The cost function is given by
\begin{equation}
J_\text{LMST} = \left(y_n - \hat y_{\text{Ten},n}  \right)^2 =:\\ e_{n}^2.
\end{equation}
The forward path and therefore the interpolated output of the tensor is straightforward and is given by 
\begin{equation}
\hat y_{\text{Ten},n} = \sum_{l_1, \ldots, l_M = 0}^{1} \sum_{r= 1}^{R} \prod_{m = 1}^{M} \matx{A}_{m,n}(k_{m,n} + l_m,r) \check{\vect{u}}_{m,n}(l_m).
\end{equation}
where in this case $k_{m,n} =\text{disc}(\hat y_{\text{LMS},n})$ is the output of the \ac{lms} algorithm \revision{ and $\text{disc}(\cdot)$ describes the discretization step according to \eqref{eq:disc}}.  
Furthermore, the \ac{lms} forward path is straightforward and given by 
\begin{align}
\hat y_{\text{LMS},n} = \mathbf w_n^{\text T} \mathbf x_{n } = \sum_{p = 1}^{P} w_{n,p} x_{n-p+1}.
\end{align}
Therefore, the cost function can be summarized as 
\begin{align}
&J_\text{LMST} =\\
& \left(y_n - \sum_{l_1, \ldots, l_M = 0}^{1} \sum_{r= 1}^{R} \prod_{m = 1}^{M} \matx{A}_{m,n}(k_{m,n} + l_m,r) \check{\vect{u}}_{m,n}(l_m)  \right)^2. \nonumber
\end{align}
Starting with the backward path of the \ac{lms} we calculate the derivative with respect to the weight vector, with the chain rule we get
\begin{equation}
\frac{\partial J_\text{LMST}}{\partial \mathbf{w}_n} = \sum_{m'=1}^M \frac{\partial J_\text{LMST}}{\partial \check{\vect{u}}_{m',n}} \frac{\partial \check{\vect{u}}_{m',n}}  {\partial \hat y_{\text{LMS},n-m'+1}  }\frac{\partial{\hat y_{\text{LMS},n-m'+1}  }}{\partial \mathbf{w}_n}.
\end{equation}
Looking at the terms individually, we calculate the derivative with respect to $\check{\vect{u}}_{m',n} $ as
\begin{align}
\frac{\partial J_\text{LMST}}{\partial \check{\vect{u}}_{m',n} } = & -2 \, e_{n} \,  \sum_{r= 1}^{R} 
\matx{A}_{m',n}(k_{m',n}:k_{m',n} +1,r) \nonumber\\
&\sum_{l_1, \ldots, l_{m'-1}, l_{m'+1}, \ldots, l_{M} = 0}^{1} \prod_{\substack{m = 1 \\ m \neq m'}}^{M} \nonumber \\
& \matx{A}_{m,n}(k_{m,n} + l_m, r) \check{\vect{u}}_{m,n}(l_m) 
\end{align}
followed by the derivative with respect to $\hat y_{\text{LMS},n-m+1} $
\begin{equation}
\frac{\partial \check{\vect{u}}_{m',n}}  {\partial \hat y_{\text{LMS},n-m'+1}  } = \begin{pmatrix} -\frac{1}{\Delta x}\\\frac{1}{\Delta x}
\end{pmatrix},
\end{equation}
and last but not least, with respect to the weights
\begin{equation}
\frac{\partial{\hat y_{\text{LMS},n-m'+1}  }}{\partial \mathbf{w}_n} = \mathbf x_{n-m'+1}^{\text{T}}.
\end{equation}
Summarizing these derivatives we get 
\begin{align}
&\frac{\partial J_\text{LMST}}{\partial \mathbf{w}_n} =   -2 \, e_{n}  \sum_{m'=1}^M    \Delta a_{m',n} \sum_{r= 1}^{R}  \mathbf v_{m',n}(r)
\mathbf x_{n-m'+1}^{\text{T}} \nonumber \\
\label{eq:derivLMST_intpol}
\end{align}
with 
\begin{align}
&\Delta a_{m',n} = \\
 &\sum_{r= 1}^{R} 
	\left(
	\frac{1}{\Delta x}  \matx{A}_{m',n}(k_{m',n} +1,r)  -\frac{1}{\Delta x} \matx{A}_{m',n}(k_{m',n} ,r) 
	\right ) \nonumber
\end{align}
and 
\begin{align}
\mathbf v_{m',n}(r)  = & \sum_{l_1, \ldots, l_{m'-1}, l_{m'+1}, \ldots, l_{M} = 0}^{1}  \prod_{\substack{m = 1 \\ m \neq m'}}^{M}  \nonumber \\
& \matx{A}_{m,n}(k_{m,n} + l_m, r) \check{\vect{u}}_{m,n}(l_m).
\end{align}
For the backward path, the derivative of $ J_\text{LMST}$ with respect to  $\matx{A}_{m',n}(k_{m',n} + l'_{m'}, :)$ is required, which can be derived as
\begin{align}
\label{eq:deriv_A_LMST_intpol}
&\frac{\partial J_\text{LMST}}{\partial \matx{A}_{m',n}(k_{m',n} + l'_{m'}, :)}  =  -2 \, e_{n} \, \check{\vect{u}}_{m',n}(l'_{m'}) \mathbf v_{m',n}
\end{align}
for all $ m' \in \{1,M\} $ and $l_{m'} \in \{0,1\}$. 
Finally, with \eqref{eq:derivLMST_intpol} and \eqref{eq:deriv_A_LMST_intpol}, we can calculate the updates of the filter coefficients and the decomposed tensor as
\begin{align}
\mathbf{w}_{n+1} =&  \mathbf{w}_{n} + 2 \mu_1  \, e_{n} \nonumber  \\
& \sum_{m'=1}^M    \,  \Delta a_{m',n} 
\sum_{r= 1}^{R} \mathbf v_{m',n}(r) \mathbf x_{n-m'+1}^{\text{T}}  
\label{eq:lmstupdate_interpolw}
\end{align}
and 
\begin{align}
\label{eq:updateA}
&\matx{A}_{m',n+1}(k_{m',n+1} + l'_{m'}, :) = \nonumber \\
& \matx{A}_{m',n}(k_{m',n} + l'_{m'}, :) +  2 \mu_2\, e_{n} \, \check{\vect{u}}_{m',n}(l'_{m'}) \mathbf v_{m',n}.
\end{align}
Note that here \eqref{eq:updateA} has to be updated for each $m' \in \{1,\dots,M\}$. 

\section{Step Size Normalization}
\label{sec:step_size}

For the purpose of improving the stability of the proposed methods, we derive a normalization for the update parameter $\mu$ for the interpolated tensor-LMS and LMS-tensor from the previous subsection. For that we consider the a priori error $e_{n+1}$ being approximated by the first order Taylor series (see~\cite{Mandic2003}).
\subsection{Tensor-LMS Normalization}


In this subsection we investigate the normalization of the tensor-LMS. Hence, we approximate the new error by the first order Taylor expansion  
\begin{align}
\label{eq:taylor_exp_intpol}
e_{n+1} \approx e_n + \sum_{r = 1}^{R}\frac{\partial e_n}{\partial \mathbf A_{m'}(:,r)}\Delta \mathbf A_{m'}(:,r)
\end{align}
with
\begin{align}
e_n = y_n - \mathbf w_n^{\text T } \hat{\mathbf {y}}_{\text{Ten},n}.
\end{align}
From \eqref{eq:updateTLMS_intpolA} it follows that
\begin{align} \label{deltaA_intpol}
&\Delta  \mathbf A_{m'}(:,r) = 2 \mu_2 e_n \mathbf S_{m',n}(:,r)
\end{align} 
and
\begin{align}
\label{eq.delta_e_delta_A}
\frac{\partial  e_n}{\partial \mathbf A_{m'}(:,r)} = \mathbf S_{m',n}^{\text T}(:,r).
\end{align}
From \eqref{eq:taylor_exp_intpol}, \eqref{deltaA_intpol} and \eqref{eq.delta_e_delta_A}  we obtain
\begin{align} \nonumber
e_{n+1} & \approx e_n  - 2 \mu_2 e_n \sum_{r = 1}^R  \mathbf S_{m',n}^{\text T}(:,r)
\mathbf S_{m',n}(:,r)\\
&= e_n  - 2 \mu_2 e_n\sum_{r = 1}^R \sum_{i = 1}^{I_{m'}}  |\mathbf S_{m',n}(i,r)|^2 \nonumber \\
& = e_n  - 2 \mu_2 e_n \norm{\mathbf S_{m',n}}_{\text F}^2,
\end{align}
\revision{
where $\norm{\cdot}_{\text F}$ denotes the Frobenius norm, }
and therefore, 
\begin{align} \label{eq:error_est_intpol}
e_{n+1} \approx  \left (1 - 2 \mu_2  \norm{\mathbf S_{m',n}}_{\text F}^2 \right )e_n .
\end{align}
Based on the Taylor expansion, 
the adaptation performance of the algorithm can be improved if the step-size is chosen to fulfill the following condition \cite{Mandic2003}
\begin{align}
\left |1 - 2 \mu_2  \norm{\mathbf S_{m',n}}_{\text F}^2 \right | < 1.
\end{align}
Hence, we obtain the following bound for the step size $\mu_2$
\begin{align}
0 < \mu_2 < \frac{1}{\norm{\mathbf S_{m',n}}_{\text F}^{2}}.
\end{align}
Be aware, that this bound has to be fulfilled for each $m' \in \{1,\dots,M\}$.
This bound can be implemented by replacing $\mu_2$ in \eqref{eq:updateTLMS_intpolA} with the time-dependent step size
\begin{align}
\mu_{2,n} = \frac{\mu_2}{\delta + \norm{\mathbf S_{m',n}}_{\text F}^{2}}, \qquad 0 < \mu_2 < 1.
\end{align}
The small regularization value ${\delta > 0}$ helps to avoid numerical issues in case of small values of the Frobenius norm.

\subsection{LMS-tensor normalization}

Similar to the tensor-LMS case, we develop the Taylor expansion for the LMS-tensor case
\begin{align}
e_{n+1} \approx e_{n} + \frac{\partial e_{n}}{\partial \mathbf{ w}_n}\Delta \mathbf{w}_n^{\text{T}}.
\end{align}
This can be done in the same way as deriving \eqref{eq:lmstupdate_interpolw} with
\begin{align}
\frac{\partial e_{n}}{\partial \mathbf{ w}_n} &=  - \sum_{m'=1}^M   \, \Delta a_{m',n} \, \sum_{r= 1}^{R}  \mathbf v_{m',n}(r)
\mathbf x_{n-m'+1}^{\text{T}} \,,
\end{align}
and
\begin{equation}
\Delta\mathbf{w}_n^\text{T} =  2 \mu_1  \, e_{n}  \sum_{m'=1}^M   \, \Delta a_{m',n}\\
\, \sum_{r= 1}^{R}  \mathbf v_{m',n}(r) \mathbf x_{n-m'+1} .
\end{equation}
Hence, we get
\begin{align}
e_{n+1} & \approx e_n  -  2 \mu_1  \, e_{n} \overbrace{\sum_{m'=1}^M   \, \Delta a_{m',n} \, \sum_{r= 1}^{R}  \mathbf v_{m',n}(r)
	\mathbf x_{n-m'+1}^{\text{T}} }^{=\mathbf{d}_n^\text{T}}\nonumber\\
&\qquad\underbrace{\sum_{m'=1}^M   \, \Delta a_{m',n} \, \sum_{r= 1}^{R}  \mathbf v_{m',n}(r)\mathbf x_{n-m'+1}}_{=\mathbf{d}_n} \nonumber \\
& = \left (1 - 2 \mu_1 \norm{\mathbf{d}_n}^2_2 \right )e_{n}\,.
\end{align}
By the same reasoning as before, the step size $\mu_1$ should be bounded by
\begin{align}
0 < \mu_1 <   \frac{1}{\norm{\mathbf{d}_n}^{2}_2}.
\end{align}
This bound can be implemented by replacing $\mu_1$ in \eqref{eq:lmstupdate_interpolw} by the time dependent step size
\begin{align}
\mu_{1,n} = \frac{\mu_1}{\delta + \norm{\mathbf{d}_n}^{2}_2}, \qquad 0 < \mu_1 < 1,
\end{align}
where $\delta$ is a small regularization value. 

\section{Computational Complexity}
\label{complexity}

\begin{table*}
	\centering
	\caption{Number of arithmetic operations per sample depending on the parameters of the considered algorithms for the estimation (forward) and update (backward) steps. \revision{In this paper we investigated the already existing classical tensor-only, TLMS, and LMST approach and the newly derived interpolated tensor-only, interpolated TLMS, and interpolated LMST.} The parameters $P$, $R$, $M$, and $I_m$ denote the length of the linear filter, rank, dimensonality of the tensors and the discretization size.}
	\resizebox{\textwidth}{!}{%
		\begin{tabular}{c|c|ccc}
			\multicolumn{2}{c|}{Algorithm}             			& Mult. & Add. & Div. \\ \hline
			
			\multirow{2}{*}{LMS}  & Forward Path                   &  $P$  & $P-1$  &  --  \\
			& Backward Path 									&  $2P+1$  & $2P$  &  $1$ \\ \hline
			
			\multirow{2}{*}{Tensor-Only}  & Forward Path   &  $(M-1)R$  & $R-1$  &  --  \\
			& Backward Path 									&  $MR(1+R(M-1)+I_m)$  & $MR(1+ I_m)$  &   $M$   \\ \hline
			
			\multirow{2}{*}{intpol Tensor-Only}  & Forward Path   &  $2^M MR$  & $2^M MR-1$  &  --  \\
			& Backward Path 									&  $2M+RM\left(1+2^{M-1}(2M-3)+I_m\right)$  & $M R (2^{M-1}+I_m)$  &  $M$  \\ \hline
			
			\multirow{2}{*}{TLMS}       & Forward Path  	    &  $P+R(M-1)$  & $P+R-2$  &  --  \\
			& Backward Path                                    &  $MR\left(P(M-1)+I_m\right)+2P(M+1)+1$  & $2P+MRI_m(P+1)$  & $1+M$  \\ \hline
			
			\multirow{2}{*}{LMST}       & Forward Path        &  $P+R(M-1)$  & $P+R-2$  &  --  \\
			& Backward Path                                     &  $1+2P+M\left(P+1+R(2M-2+I_m)\right)$  & $M\left(R(3+I_m)+P-1\right)+P$  & $1+M$  \\ \hline
			
			\multirow{2}{*}{intpol TLMS}       & Forward Path  	    &   $2^MRM+P$   & $2^MRM +P-2$  &  --  \\
			& Backward Path                                    &  $MRP \left(3+2^{M-1}(2M-3)\right)+I_mR(M+1)+M^2+2P$  & $I_m R(1+M)+P (2^{2M-1}+2^{M-1}(1-R)+2) -M$  & $1+M$ \\ \hline
			
			\multirow{2}{*}{intpol LMST}       & Forward Path        &  $2^MRM+P$   & $2^MRM +P-2$  &  --  \\
			& Backward Path                                     &  $P(M+2) + 2(M+1) +MR \left(2 + 2^{M}(2M-3)+ I_m\right) $  & $2P+MR( 2^{M-1}3+I_m+1)-1$  & $1+M$  
	\end{tabular}}
	\label{tab:compl1}
\end{table*}

In this section we compare the complexity of all methods covered in this paper in terms of additions, multiplications and divisions. For the \ac{lms} algorithm one can see in \tabref{tab:compl1}, that the update requires $2P+1$ multiplications and $2P$ additions. For the numbers of required operations in the forward and backward paths of the classical tensor the equations in Section~\ref{sec:basics} are investigated. For our proposed interpolated tensor models a close look at the equations derived in Section~\ref{sec:sys_ident} and also on the step size normalization in \secref{sec:step_size} is taken. For simplicity, we chose for each dimension the same discretization size, hence
\begin{align}
	l_0 = l_1 =  \dots = l_{m}.
\end{align} 
For the forward path of the interpolated tensor-only approach according to~\eqref{eq:estimate_intpol} $(M-1)R$ multiplication and $R-1$ additions are needed. 
Furthermore, the backward path for the interpolated SGD update given in~\eqref{eq:update_intpol} is more complex and summarized in Table~\ref{tab:compl1}. Moreover, the complexity for the interpolated TLMS and LMST algorithms is derived in a similar way by looking at the corresponding equations in Section~\ref{sec:sys_ident} and can also be found in Table~\ref{tab:compl1} in a general form.
Evidently, the rank $R$ and the dimension $M$ of both, the classical tensor approaches as well as the interpolated methods affect the required number of operations per sample the most. 
While the computational complexity of classical tensor methods grows linearly with the dimension $M$, the computational complexity of the interpolated tensor methods grows exponentially with $M$. 

\section{Performance Simulations}
\label{sec:results}

\begin{table*}
	\centering
	\caption{\revision{For the six different experiments investigated, each of the applied algorithms needed different parameters. These parameter values of the classical tensor methods and the interpolated tensor methods used in the experiments are summarized here.}}
	\resizebox{\textwidth}{!}{%
		\begin{tabular}{c|cccc|cc|cccc|cc|cccc|cccc}
			\multirow{2}{*}{Experiment} & \multicolumn{4}{c|}{Tensor} & \multicolumn{2}{c|}{LMS} & \multicolumn{4}{c|}{intpol Tensor } & \multicolumn{2}{c|}{intpol LMS} & \multicolumn{4}{c|}{Tensor-Only} & \multicolumn{4}{c}{intpol Tensor-Only}                                \\
			& $\mu_2$ &  $R$ &       $M$       & $I$ &  $\mu_1$ &   $P$   & $\mu_2$ & $R$ &       $M$       & $I$ & $\mu_1$ &   $P$   & $\mu$ & $R$ &        $M$         & $I$ & $\mu$ & $R$ &        $M$         &  $I$ \\ \hline
			1      &  0.009  &  1  &        1       &   50   &  0.009  &         7      &   0.01   &  1  &        1         &   10     &   0.01   &         7          & 0.05  & 50  &         7    &   25  &   0.1    &   10   &  3  &        10                \\
			2      &  0.1  &  1  &        1         &   50         &  0.0075  &         5    &   0.01 &   1   &  1  &        10       &   0.001   &         5          & 0.01  & 50  &         5  &   23  &  0.1  &        10        &   3  &   10         \\
			3      &  0.008  &  1  &        3        &   50          &  0.005  &         5     &   0.01  &   10 &  2  &        16       &   0.01  &         5          & 0.025  & 100  &         7&   25   &  0.4 &        20            &   3  &   10           \\
			4      &  0.05  &  10  &        2        &   16          &  0.002  &         7   &  0.01   &   10  &   2 &        16        &   0.01   &         7         & 0.05  & 200  &         8    &   25  &  0.1  &        20       &   3  &   32  \\
			5      &  0.1  &  30  &        3        &   50          &  0.02  &         3   &   0.1  &   16 &  3  &        20         &    0.8  &         3          & 0.09  & 40  &         3&   30   &  0.1 &        40            &   3  &   20           \\
			6      &  0.4  &  10  &        4        &   40          &  0.001  &         4  &  0.5   &   16  &   2 &        30         &   0.001   &         2         & 0.4  & 20  &         3    &   20  &  0.4  &        20       &   3  &   10  
1	\end{tabular}}
	\label{tab:simsettings}
\end{table*}

We conclude our investigations on interpolated tensors with performance simulations by comparing them with the classical tensor approach proposed in~\cite{tensor}. The general aim of all considered algorithms is to minimize the error between the wanted signal $d_{n}^{(l)}$ and its replica $\hat{y}_{n}^{(l)}$, where $n$ is the time and $l$ indicates the run (i.e., $l$-th  repetition). Therefore, the metric we use for comparison is the \ac{nmse} 
\begin{equation}
\text{NMSE}_\text{dB} = 10\log_{10}\left(\frac{1}{L}\sum_{l=1}^{L}\frac{\left(d_{n}^{(l)}-\hat{y}_{n}^{(l)}\right)^2}{\left(d_{n}^{(l)}\right)^2}\right)\,,
\end{equation}
where $L$ is the total number of runs of a given experiment. 
Note that in $d_n$ the measurement noise is not included cf. \figref{fig:overviewarch}.

In order to validate the proposed interpolated tensor architecture, several experiments were performed. 
 The experiments 1-4 conducted in this work formally correspond to those conducted in \cite{tensor,scrpiniti_cascade} . The description of these experiments is repeated in the following for the sake of completeness. 
All previously discussed methods, the interpolated tensor only, the interpolated TLMS and the interpolated LMST are evaluated in their performance. 
The settings (step size, rank, etc.)  have been chosen empirically, such that the all algorithms perform the best.
For all experiments the weights and the \ac{tdl} of all adaptive algorithms are initialized with zeros and the tensors are initialized with standard normal distributed random values with variance $0.01$.  For the experiments a  \ac{snr} of $\SI{10}{\dB}$ was used. In Example 1 and 2 we want to stress the adaptive nature of the considered algorithms. Therefore, after half of the simulation time we change the linear part of the unknown system. A reason for this happening can be temperature drifts within the system or simply because of switching to other linear filters during operations.

\subsection{Experiment 1}
The first experiment focuses on the identification of a Hammerstein model, i.e., the system includes a static nonlinearity followed by a linear part.  
\revision{The input signal $x_n$ 
is generated by the
relationship }
\begin{equation}
x_n = ax_{n-1}+\sqrt{1-a^2}\nu_n\,,
\label{eq:modelx}
\end{equation}
where $\nu_n$ is a zero mean white Gaussian noise with unitary variance and $0\leq a<1$ where the parameter $a$ determines the level of correlation between adjacent samples, for further details on the choice of $a$ we refer to~\cite{tensor,scrpiniti_cascade}.
The first block of the system is described by a static nonlinearity in the form of 
\begin{equation}
y_n^\text{nl} = \frac{2x_n}{1+\abs{x_n}^2}\,.
\end{equation}
This part of the system imitates the saturation behavior of a \ac{pa} in a satellite communication ~\cite{720369}, cf. the first experiment in\cite{tensor}.
The second block is an unknown linear \ac{fir} filter  $\mathbf{w}_0\in\mathbb{R}^7$ which models a part of the transmission path.
In this experiment we change the coefficient of this linear system after half of the run time. 
The $\mathbf{w}_0$ coefficients are initially given by 

\revision{
$\begin{pmatrix}
 0.9 & -0.6 & 0.3 & -0.15 & 0.1 & -0.025 & 0.005
 \end{pmatrix} ^{\text{T}}$} as in~\cite{scrpiniti_cascade} and later by \revision{$
 \begin{pmatrix}
 0.6  -0.4  0.25  -0.15  0.1  -0.05  0.001
 \end{pmatrix}^{\text{T}}$}, cf.~\cite{tensor}.

The interpolated tensor-only approach and, because of the Hammerstein structure, the interpolated TLMS architecture are compared to the classical tensor-only and  TLMS  approach. 
In Fig.~\ref{fig:sim1} the performance in terms of \ac{nmse} for all four algorithms is depicted. One can see that both proposed interpolated variants outperform their classical counterparts significantly. Even though both TLMS versions have rank and dimensionality equal to one the complexity is lower in the interpolated case because a lower discretization is sufficient (cf. Table~\ref{tab:simsettings} and~\ref{tab:compl2}). For the tensor-only approach the classical method needs way more dimensionality and rank than its interpolated counterpart. Nevertheless, almost all approaches react quickly and as expected to the change of the linear subsystem mentioned previously, compare to~\cite{tensor}. The exception is the classical tensor-only approach which not only requires by far the most arithmetic operations but also performs worse than the other approaches. This approach is neither performing well enough nor adapting to the change of the system appropriately. It performs even worse after the change. 

\subsection{Experiment 2}
The second experiment follows a Wiener model. Here, the general setup consists of the same nonlinearity as in the previous experiment with the difference that the nonlinearity occurs after the linear filter consisting of only 5 taps, for further details we refer to ~\cite{tensor}.   
Hence, for this example the LMST structure is used. The performance of the algorithms is displayed in Fig.~\ref{fig:sim2}. As can be seen, again both proposed methods outperform their classical counterparts. 
Similar to the previous experiment, the rank $R$ and dimensionality $M$ of the tensor in the LMST variants are set to one, while the classical tensor-only approach requires way higher order and rank than the interpolated version. As in the previous example all algorithms besides the classical tensor-only approach adapt to the change of the system as expected and yield similar performance after the change.  

\subsection{Experiment 3}
The third experiment again has a Wiener model setup with a 5 tap linear filter, but we choose a different nonlinearity  
\begin{equation}
y^\text{nl}_n = \sin^2\left(x_n\right) + \sin^3\left(x_{n-1}\right)+\sin^4\left(x_{n-2}\right)\,.
\end{equation}
Instead of the simple one in experiment two, the output of the nonlinearity at time $n$ depends on the current and last two input values of this block. Hence, the nonlinearity has a memory. The input signal $x_n$ follows the same relationship as in Experiment 1 in  \eqref{eq:modelx}.
Also for this example,  Fig.~\ref{fig:sim5} shows that the proposed method outperforms the classical tensor solutions. The interpolated LMST performs the best followed by the interpolated tensor-only solution. Complexity wise the interpolated LMST is more expensive then its classical counter part. In comparison the interpolated tensor-only solution needs less complexity then the classical tensor-only variant (cf. Table~\ref{tab:compl2}).

\subsection{Experiment 4}
The fourth experiment follows again a Hammerstein model with the same setup as experiment 3 with the linear part, consisting of $7$ taps and nonlinear part swapped. The evaluation in Fig.~\ref{fig:sim6} shows that the interpolated tensor-only outperforms the classical tensor-only solution significantly in terms of \ac{nmse} and convergence speed (by means of the slope of the \ac{nmse} curves at the beginning of the adaptation) while the interpolated combined version is only slightly better than the classical one. Regarding complexity the LMST algorithms are only slightly different. While the classical method needs more additions, the interpolated shows a higher number of multiplications.    
 
 \subsection{Experiment 5}
 The fifth experiment is quite similar to the first one and has also a Hammerstein structure. Beside that, the system is motivated by a specific mobile communication problem instead of satellite communications. The nonlinearity is inspired by the Tx harmonics problem \cite{Auer2021} and also imitates the nonideality of a \ac{pa}. It is modeled by
 \begin{equation}
 y_n^\text{nl} = x_n^2\, ,
 \end{equation}
 followed by a linear filter of length $3$.
 The input signal does not follow the relation \eqref{eq:modelx} but is white Gaussian noise with variance one.
 The results of this simulation can be seen in Fig.~\ref{fig:sim7}. The proposed approaches significantly outperform the classical tensor methods. While the interpolated TLMS variant performs slightly worse than the interpolated tensor-only solution it's advantage is that it need less arithmetical operations.  
  
 \subsection{Experiment 6}
\begin{figure*}
	\centering
	\begin{subfigure}[t]{0.39\textwidth}
		\centering
		\includegraphics[width=\textwidth]{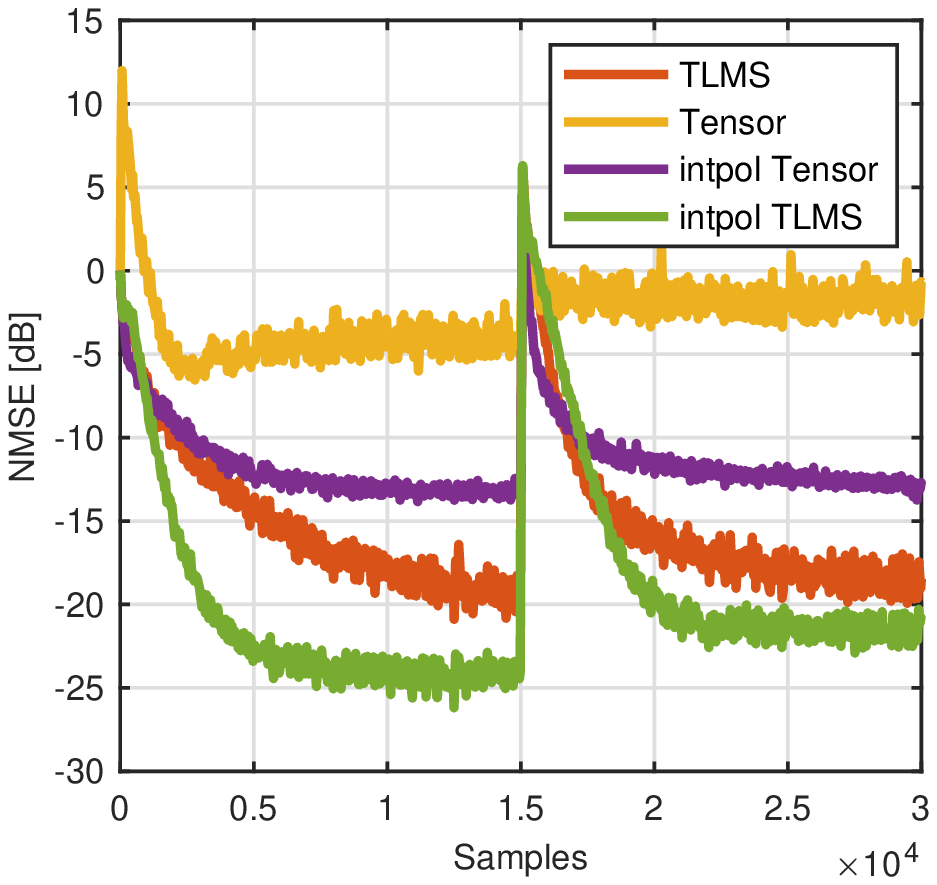}
		\caption{}
		\label{fig:sim1}
	\end{subfigure}%
	~ 
	\begin{subfigure}[t]{0.39\textwidth}
		\centering
		\includegraphics[width=\textwidth]{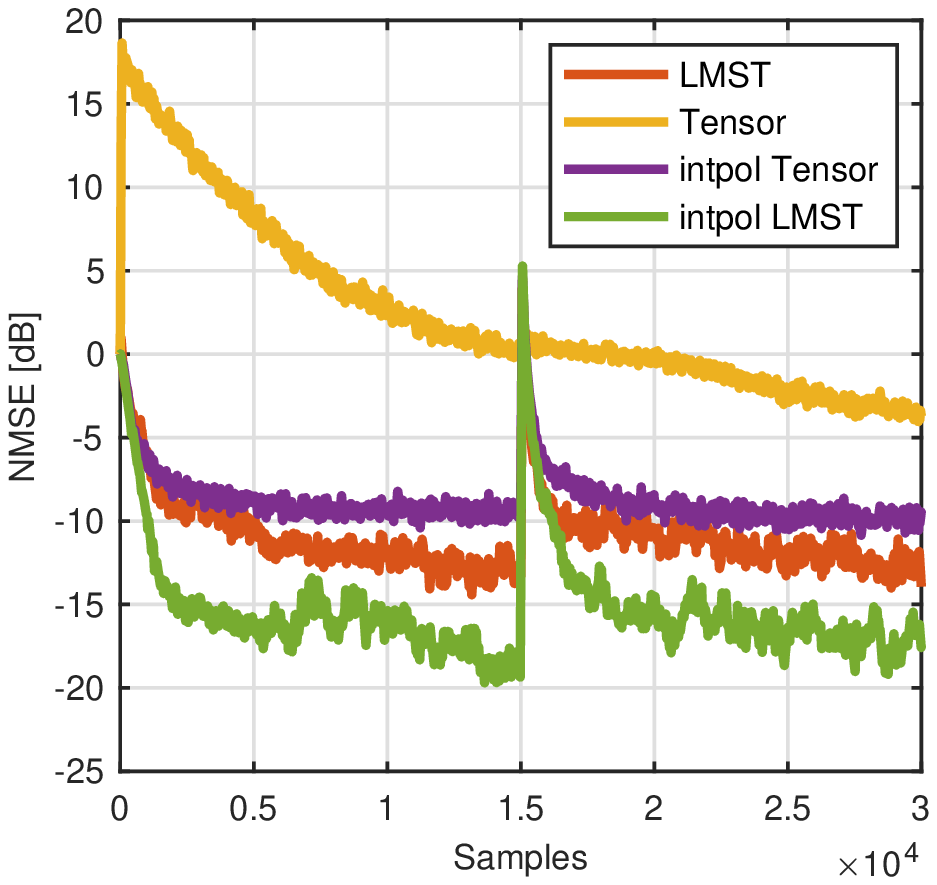}
		\caption{}
		\label{fig:sim2}
	\end{subfigure}
	~
	\begin{subfigure}[t]{0.39\textwidth}
		\centering
		\includegraphics[width=\textwidth]{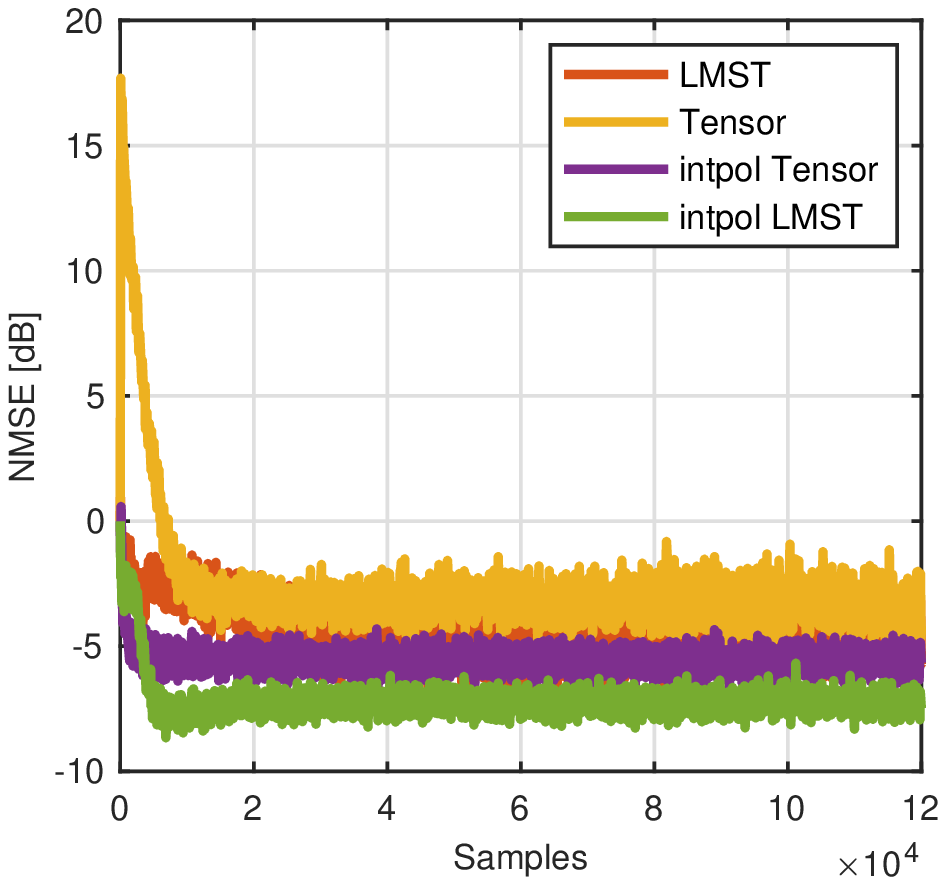}
		\caption{}
		\label{fig:sim5}
	\end{subfigure}%
	~
	\begin{subfigure}[t]{0.39\textwidth}
		\centering
		\includegraphics[width=\textwidth]{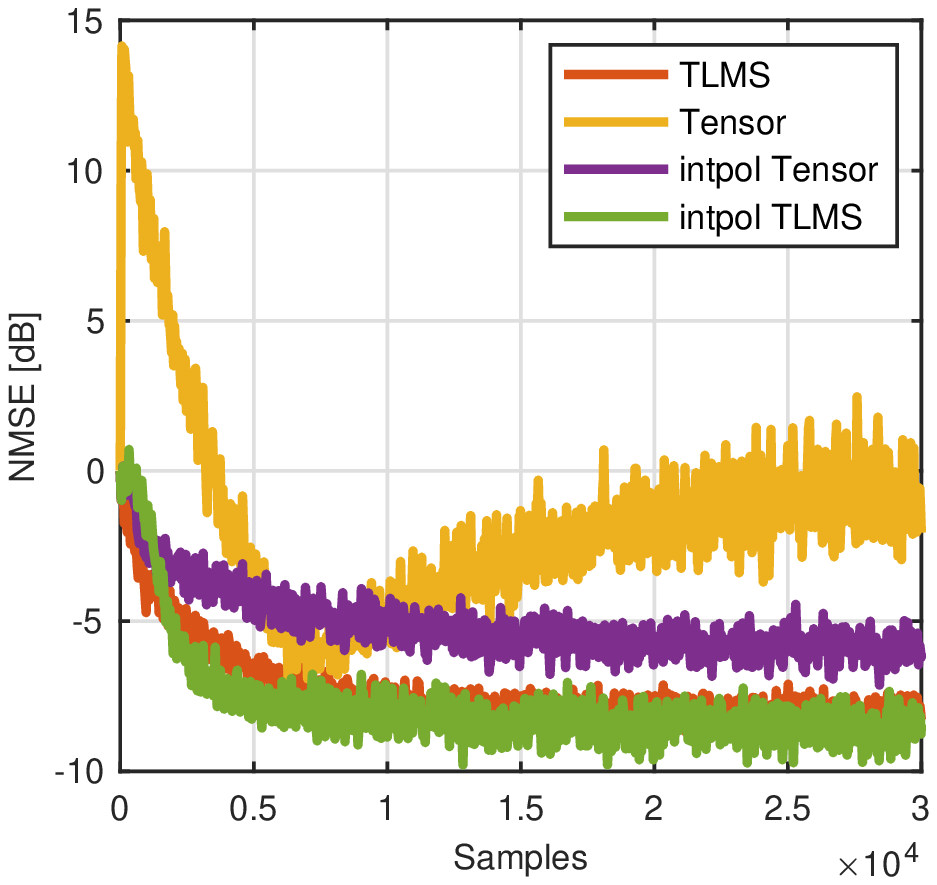}
		\caption{}
		\label{fig:sim6}
	\end{subfigure}
	~
	\begin{subfigure}[t]{0.39\textwidth}
		\centering
		\includegraphics[width=\textwidth]{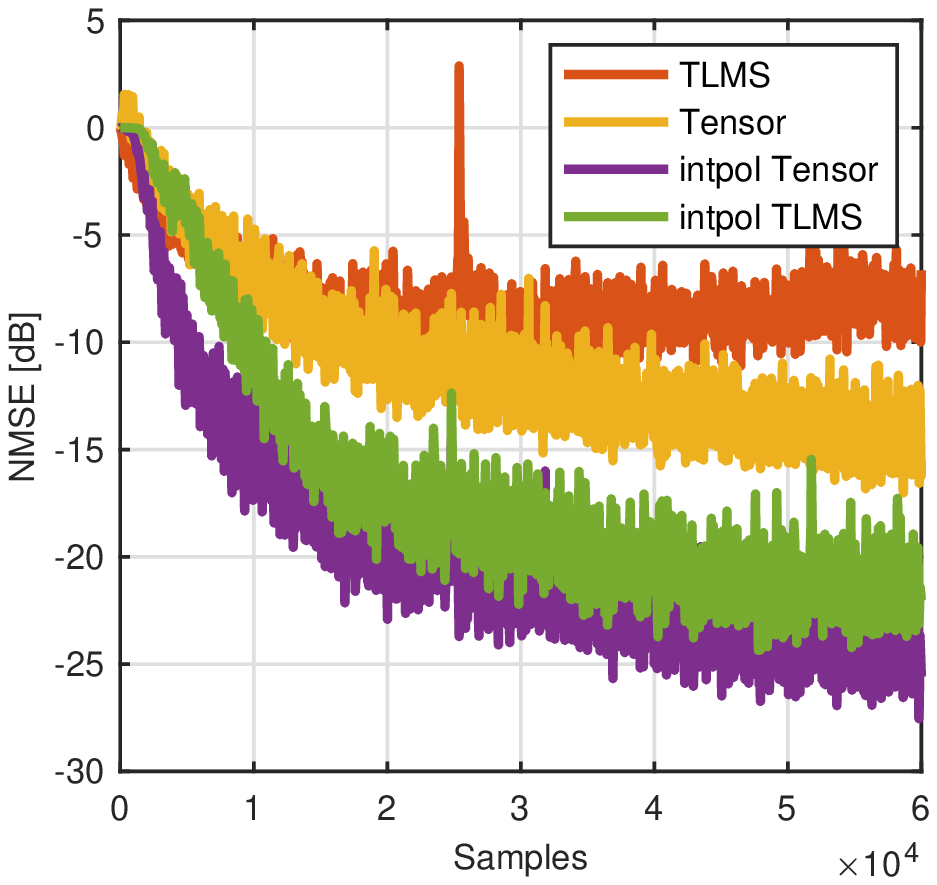}
		\caption{}
		\label{fig:sim7}
	\end{subfigure}
	~
	\begin{subfigure}[t]{0.39\textwidth}
		\centering
		\includegraphics[width=\textwidth]{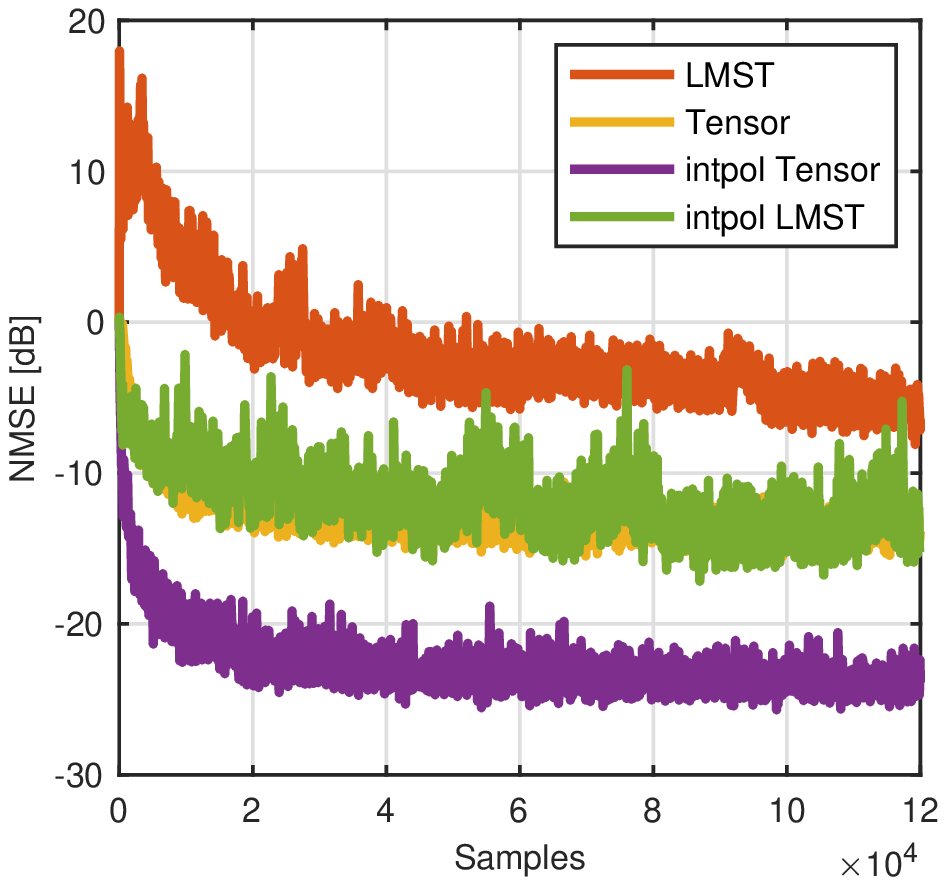}
		\caption{}
		\label{fig:sim8}
	\end{subfigure}
	\caption{Comparison NMSE convergence for all considered algorithms for (a) the Hammerstein case, (b) the Wiener case, (c)  the Wiener case with a nonlinearity containing memory, and (d) the Hammerstein case with a nonlinearity depending on past samples, (e) the Hammerstein case inspired by mobile communications, (f) the Wiener model inspired by mobile communications}
	\label{fig:sims}
\end{figure*}

The sixth experiment is comparable to the second experiment following a Wiener structure. It is again motivated by a specific mobile communication problem, namely the \ac{imd2} problem \cite{Auer2020}. The linear filter has a $3$ taps and the nonlinearity has the following form 
 \begin{equation}
 y_n^\text{nl} = |x_n|^2\,.
 \end{equation}
The input data is again chosen to be white Gaussian noise with variance one. The evaluation in Fig.~\ref{fig:sim8} shows that both interpolated methods outperform their classical counterparts and all solutions have a fast convergence speed. 

\subsection{Comparison}

%
%
%

\begin{table*}
	\centering
	\caption{Complexity in terms of multiplications, additions and divisions of the considered algorithms in all experiments conducted.}
	\begin{tabular}{c|c|ccccc}
		Experiment & Operation            & Tensor-only & TLMS/LMST & intpol Tensor-only & intpol TLMS/LMST \\ \hline
		\multirow{3}{*}{1}  &   Mult.    &  114400  & 86 &  936  & 58 \\
		& Add.       &  9149  & 420 &  659  & 54 \\
		&  Div.        &  7  & 2  &  3  & 2 \\ \hline
		
		\multirow{3}{*}{2}  &   Mult.    &  56200  & 72  &  1866  & 36  \\
		& Add.       &  6049  & 66 &  1319  & 28  \\
		&  Div.        &  5  & 2  &  3  & 2\\ \hline
		
		\multirow{3}{*}{3}  &   Mult.    &  438800  & 198  &  1866  & 551 \\
		& Add.       &  18299  & 180  &  1319  & 512 \\
		&  Div.        &  7  & 4  &  3  & 3 \\ \hline
		
		\multirow{3}{*}{4}  &   Mult.    &  2283000  & 520 &  3186  & 1285   \\
		& Add.       &  41799  & 2589  &  2639  & 467 \\
		&  Div.        &  8  & 3  &  3  & 3  \\ \hline
		
		\multirow{3}{*}{5}  &   Mult.    &  12200  & 1308  &  4926  & 3842 \\
		& Add.       &  2559  & 3863  &  3839  & 1328 \\
		&  Div.        &  3  & 4  &  3  & 4 \\ \hline
		
		\multirow{3}{*}{6}  &   Mult.    &  3100  & 1053 &  1866  & 1296   \\
		& Add.       &  679  & 1076  &  1319  & 1251 \\
		&  Div.        &  3  & 3  &  3  & 3  
	\end{tabular}
	\label{tab:compl2}
\end{table*}

For each experiment, the number of required operations is displayed in Table~\ref{tab:compl2}. This shows that the TLMS/LMST versions are always less complex than the tensor-only counterparts. Furthermore, it can be seen that the interpolated tensor-only method requires less operations than the classic version in all simulations, while for the Wiener and Hammerstein approaches it varies from example to example. Sometimes, the interpolated case requires less operations. In all cases, the interpolated versions outperform the classical approaches in terms of the final NMSE. Only in experiment 4 the TLMS method almost performs as good as the interpolated TLMS method. An evaluation of the convergence speed is more difficult due to the fact that most of the considered methods adapt similarly. Two clear exceptions can be identified in form of the classical tensor-only approach, which suffers from convergence issues, and in form of experiment 6 for which the classical LMST method suffers from a bad adaptation behavior as well.  In general one can say that the classical TLMS/LMST variants seem to deal worse with normally distributed inputs compared to the input samples with the structure from \eqref{eq:modelx} (cf. experiment~1 - 4).
All in all, our proposed interpolated TLMS/LMST methods are low-cost architectures that can achieve better performance than the state-of-the-art its classical counter parts. 

\section{Conclusion}
\label{conclusion}
In this paper we presented three novel adaptive learning schemes which combine the state-of-the-art tensor methods with interpolation for the general problem of system identification. These three schemes are an interpolated tensor-only, approach, a Hammerstein and a Wiener system learning algorithm. 
The latter two methods are based on the well-known LMS algorithm. Several experiments demonstrate the effectiveness of the proposed approaches in terms of \ac{nmse}. Moreover, a constraint on the choice of the learning rate is also derived, in order to ensure stability. 
Additionally, a complexity analysis has been carried out and shows that the proposed interpolated tensor-only method requires less operations than the classical method for all six experiments, and even the interpolated TLMS/LMST methods show almost the same complexity.


\revision{
\bibliographystyle{IEEEtran}
\bibliography{mybib}
}



\begin{IEEEbiography}[{\includegraphics[width=1in,height=1.25in,clip,keepaspectratio]{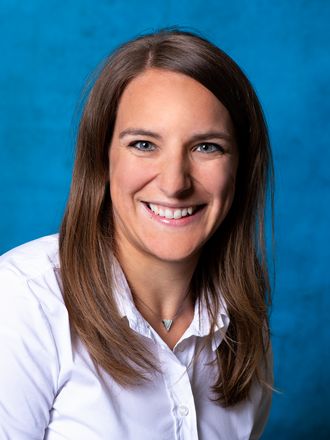}}]{Christina Auer} (GSM'18)
		was born in Linz, Austria in 1987. From 2005 till 2012 she studied Applied Mathematics at the Johannes Kepler University Linz. The topic of her master thesis, she partly wrote at the UCLA, was image denoising. After graduation, she joined the company LCM (Linz Center of Mechatronics) where she was working in the R\&D, especially with Kalman Filtering for localization. Since September 2017, she is a member of the Institute of Signal Processing at the Johannes Kepler University Linz. Now she is working towards her Ph.\,D.\ as a member of the Christian Doppler laboratory for Digitally Assisted RF Transceivers for Future Mobile Communications. Her research activities focus on self-interference cancellation using kernel adaptive filtering and machine learning.
\end{IEEEbiography}

\begin{IEEEbiography}[{\includegraphics[width=1in,height=1.25in,clip,keepaspectratio]{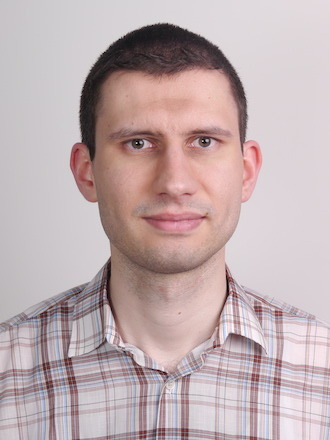}}]{Thomas Paireder} (GSM'19)
	received his bachelor’s degree (Hons.) in information electronics and his master’s degree (Hons.) in electronics and information technology from Johannes Kepler University Linz (JKU), Linz, Austria, in 2016 and 2018, respectively. In his master’s thesis, he implemented a real-time processing system for laser-ultrasonic signals, which is utilized in nondestructive testing. He is currently pursuing the Ph.\,D.\ degree at the Institute of Signal Processing, JKU. His topic is receiver interference cancellation by means of adaptive signal processing methods.
\end{IEEEbiography}

\begin{IEEEbiography}[{\includegraphics[width=1in,height=1.25in,clip,keepaspectratio]{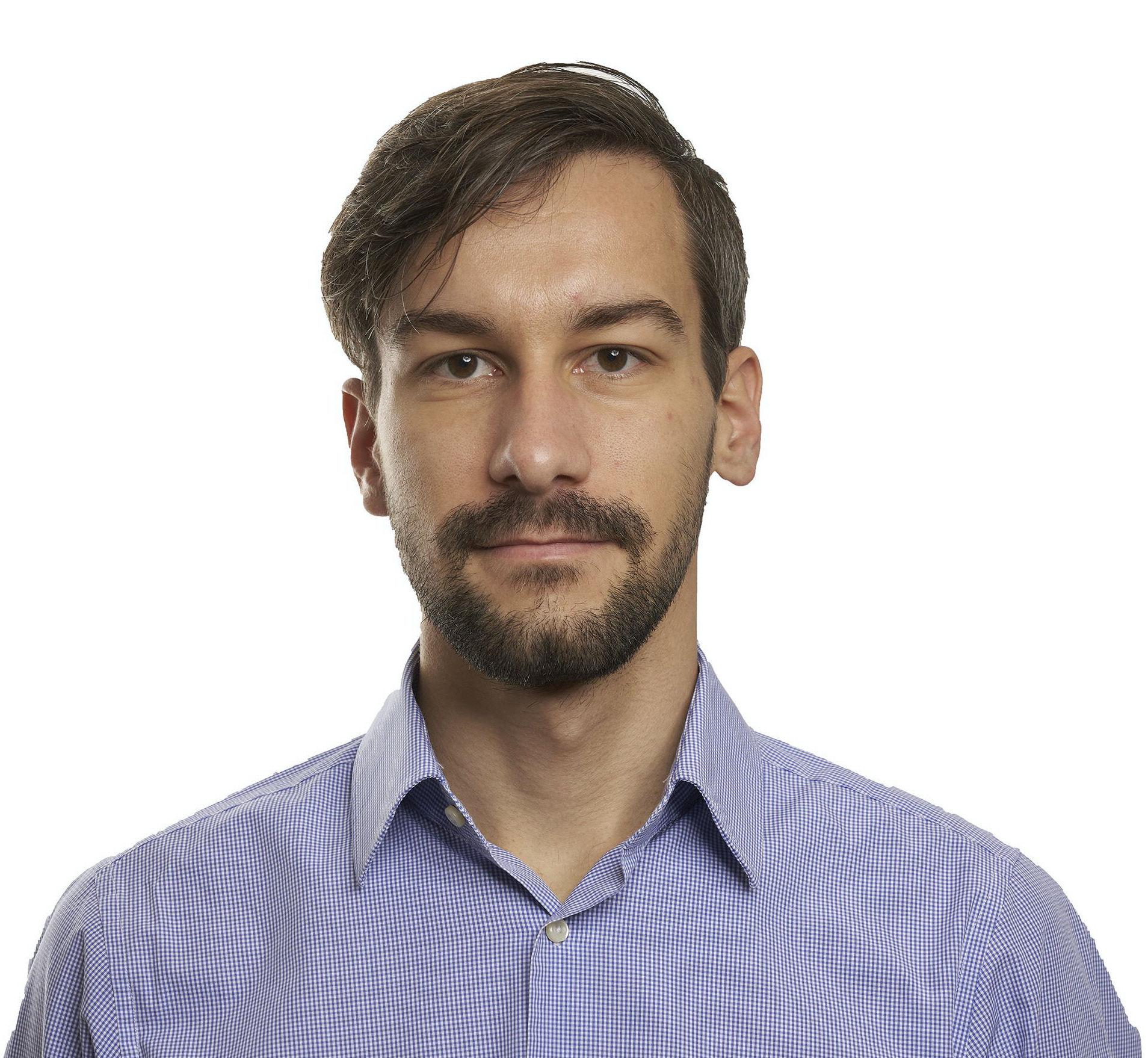}}]{Oliver Ploder} 
	(S'18) was born in Graz, Austria, in 1992. From 2012 till 2016 he studied at the Johannes Kepler University of Linz where he obtained his bachelor's degree in Information Electronics. He received his master's degree in Telecommunication engineering with a focus on Wireless Communications from the Universitat Politècnica de Catalunya (BarcelonaTech) in 2018. His master thesis was written in cooperation with the Networked and Embedded Systems Laboratory (NESL) at UCLA on the topic of secure state estimation in Cyber Physical Systems (CPS). Since September 2018, he is a member of the Institute of Signal Processing at the Johannes Kepler University, where he is working towards his PhD,  focusing his research on receiver interference cancellation for LTE and LTE-A RF transceiver systems by means of machine learning.
\end{IEEEbiography}

\begin{IEEEbiography}[{\includegraphics[width=1in,height=1.25in,clip,keepaspectratio]{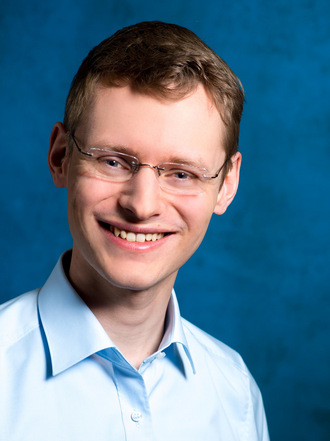}}]{Oliver Lang}
	(S’14–M’17) was born in Schärding, Austria in 1987. He studied Electrical Engineering at the Vienna University of Technology and received his bachelor's degree in 2011. In the next two years Oliver Lang studied Microelectronics at the Vienna University of Technology and finished it with excellent success. The topic of his master thesis was the development and analysis of models for a scanning microwave microscope. From 2014 until 2018, he was a university assistant at the Institute of Signal Processing, Johannes Kepler University Linz, where he received his Ph.D. in 2018. The title of his dissertation is “Knowledge-Aided Methods in Estimation Theory and Adaptive Filtering”, where he invented several interesting estimators and adaptive filters for special applications. From 2018 until 2019 he worked at DICE GmbH in Linz, which is a subsidiary company of Infineon Austria GmbH. Since March 2019, he is a university assistant with Ph.D. at the Institute of Signal Processing.
\end{IEEEbiography}

\begin{IEEEbiography}[{\includegraphics[width=1in,height=1.25in,clip,keepaspectratio]{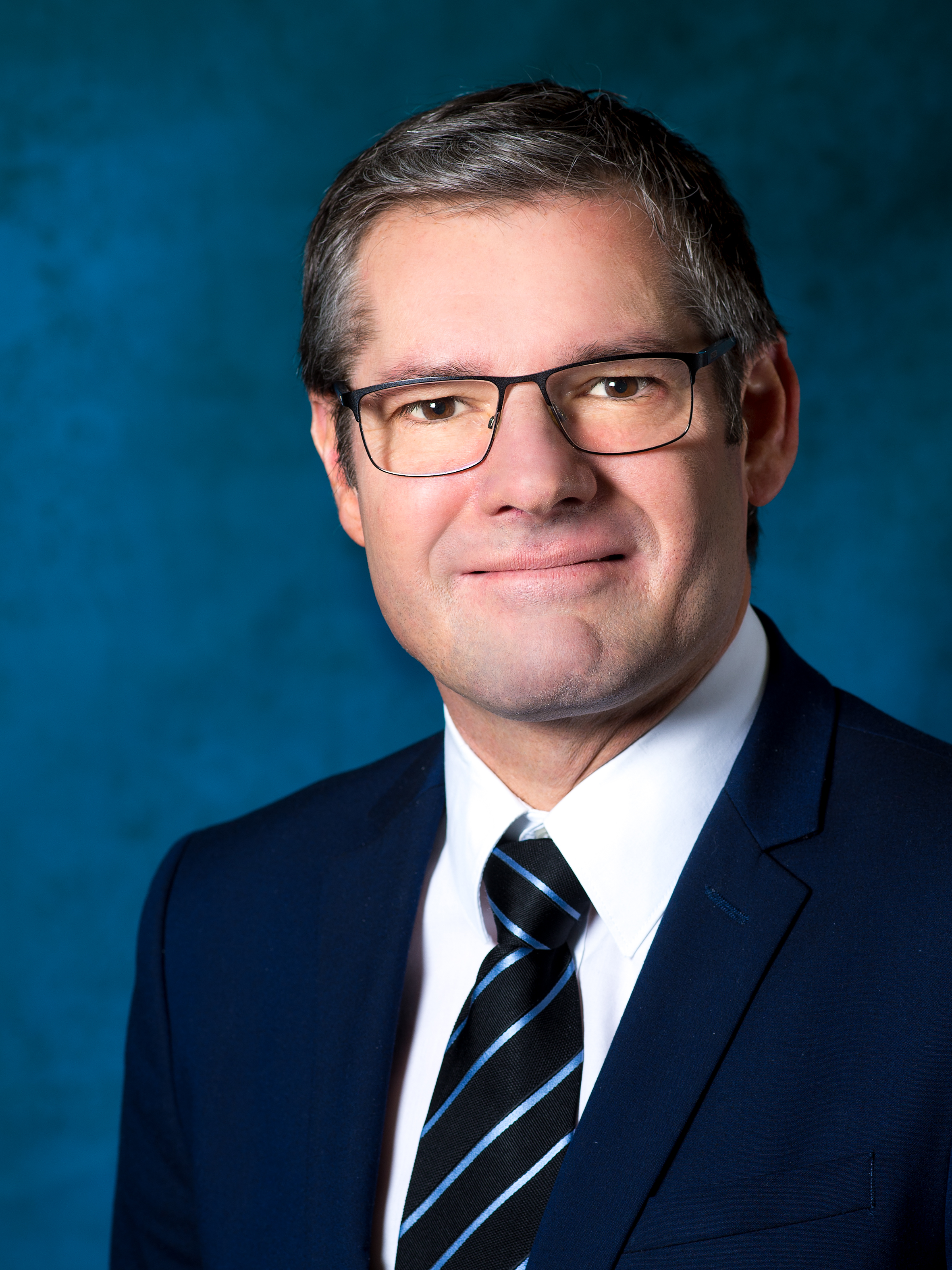}}]{Mario Huemer} (M'00--SM'07)
	received his Dipl.-Ing.\ and Dr.~techn.\ degrees from the Johannes Kepler University (JKU) Linz, Austria, in 1996 and 1999, respectively. After holding positions in industry and academia he became an associate professor at the University of Erlangen-Nuremberg, Germany, from 2004 to 2007, and a full professor at Klagenfurt University, Austria, from 2007 to 2013. Since September 2013, Mario Huemer is heading the Institute of Signal Processing at JKU Linz as a full professor, and since 2017 he is co-head of the ``Christian Doppler Laboratory for Digitally Assisted RF Transceivers for Future Mobile Communications''. His research focuses on statistical and adaptive signal processing, signal processing architectures, as well as mixed signal processing with applications in information and communications engineering, radio frequency transceivers for communications and radar, sensor and biomedical signal processing. Within these fields he has published more than 270 scientific papers. From 2009 to 2015 he was member of the editorial board of the International Journal of Electronics and Communications (AEU), and from 2017 to 2019 he served as an associate editor for the IEEE Signal Processing Letters. Mario Huemer has received the dissertation awards of the German Society of Information Technology (ITG) and the Austrian Society of Information and Communications Technology (GIT), respectively, the Austrian Kardinal Innitzer award in natural sciences, and the German ITG award.
\end{IEEEbiography}


\end{document}